\documentclass[12pt]{article}

\usepackage[T1]{fontenc} \usepackage[english]{babel}
\usepackage{latexsym, amssymb}

\begin{document}

\newtheorem{theo}{Theorem}
\newtheorem{defi}[theo]{Definition}
\newtheorem{theo1}[theo]{Theorem}
\newtheorem{prop}[theo]{Proposition}
\newtheorem{lemm}[theo]{Lemma}
\newtheorem{coro}[theo]{Corollary}

\newcommand{\R}{\mbox{${\bf R}$}}
\newcommand{\Z}{\mbox{${\bf Z}$}}
\newcommand{\Q}{\mbox{${\bf Q}$}}
\newcommand{\N}{\mbox{${\bf N}$}}
\newcommand{\elt}{\mbox{$v \in G$}}
\newcommand{\elg}{\mbox{$\phi \in \Gamma$}}
\newcommand{\gstab}{\mbox{$G_v$}}
\newcommand{\caa}{\mbox{$C_1(\alpha)$}}
\newcommand{\cia}{\mbox{$C_i(\alpha)$}}
\newcommand{\ima}{\mbox{${\rm Im}$}}
\newcommand{\de}{\mbox{${\rm deg}$}}
\newcommand{\dist}{\mbox{${\rm dist}$}}
\newcommand{\wre}{\mbox{${\rm Wr}$}}
\newcommand{\res}{\mbox{$\mid$}}
\newcommand{\nin}{\mbox{$\ \not\in\ $}}
\newcommand{\fix}{\mbox{${\rm fix}$}}
\newcommand{\pr}{\mbox{${\rm pr}$}}
\newcommand{\sym}{\mbox{${\rm Sym}\;$}}
\newcommand{\aut}{\mbox{${\rm Aut}\;$}}
\newcommand{\auttn}[0]{\mbox{${\rm Aut}\;T_n$}}
\newcommand{\autx}[0]{\mbox{${\rm Aut}\;X$}}
\newcommand{\QZ}{\mbox{${\rm QZ}$}}

\title{Normal subgroups of groups acting on trees and automorphism groups of graphs}
\author{R\"ognvaldur G. M\"oller\footnote{Science Institute,University of Iceland, IS-107 Reykjav\'ik, Iceland.  e-mail: roggi@raunvis.hi.is}\  
\mbox{ and} Jan Vonk\footnote{Mathematical Institute, 24-29 St Giles',
Oxford, OX1 3LB, England.  e-mail:  Jan.Vonk@maths.ox.ac.uk.}
\footnote{The second author wishes to thank the ERASMUS programme for support.}
}

\maketitle

\begin{abstract}
\noindent
Let $T$ be a tree and $e$ an edge in $T$.  If $C$ is a component of $T\setminus e$ and both $C$ and its complement are infinite we say that $C$ is a half-tree.  The main result of this paper is that if $G$ is a closed subgroup of the automorphism group of $T$ and $G$ leaves no non-trivial subtree invariant and fixes no end of $T$ then the subgroup generated by the pointwise stabilizers of half-trees is topologically simple.  This result is used to derive analogues of recent results of Caprace and De Medts \cite{CapraceDeMedts2011} and it is also applied in the study of the full automorphism group of a locally finite primitive graph with infinitely many ends.
\end{abstract}

\bigskip\bigskip\bigskip

Running title:  NORMAL SUBGROUPS OF GROUPS ACTING ON TREES

\newpage

\section*{Introduction}

In  the first half of this paper we study a variant of Tits' simplicity result from his ground breaking paper \cite{Tits1970} on automorphism groups of trees.   Tits studies a group action on a tree $T$ such that the action satisfies a certain independence property called  property P.  Here we look at a different property which can also be thought of as an independence property.  This property is defined in terms of {\em half-trees}.  If $e$ is an edge in a tree $T$ and both components of $T\setminus e$ are infinite then we call these components half-trees.  A group acting on $T$ has property H if the pointwise stabilizer of every half-tree is non-trivial.  For a closed subgroup $G$ of the automorphism group of a tree $T$ we let $G^{++}$ denote the closure of the subgroup generated by all pointwise stabilizers in $G$ of half-trees.
The main result of this paper is 

\medskip

{\bf Theorem~\ref{TSimple}}\ \ 
{\em Let $G$ be a closed subgroup of the automorphism group of some tree $T$.   Assume that no proper non-empty subtree of $T$ is invariant under $G$ and no end of $T$ is fixed by $G$.  Suppose also that $G$ has property H.
If $N$ is a non-trivial closed subgroup of $G$ normalized by $G^{++}$ then $N$ contains $G^{++}$.   In particular, the subgroup $G^{++}$ is topologically simple. } 

\medskip

\noindent
It is also shown that many of the results in a recent paper of Caprace and De Medts \cite{CapraceDeMedts2011} that are proved for groups satisfying property P are also true if property H is assumed instead.

In the second half of the paper the full automorphism group of a graph with infinitely many ends is studied.  Such groups satisfy a certain independence property because if $A$ is a set of vertices or edges in a graph $X$ such that $X\setminus A$ is not connected then the subgroup of the full automorphism group fixing all the elements in $A$ and leaving invariant each component of $X\setminus A$ acts on each component independently of what it does on the other components.  First we look at the automorphism group of a transitive graph with connectivity 1, i.e. connected graphs were the removal of a single vertex produces a disconnected graph.  The automorphism group of such a graph can be studied with the aid of Tits' original result and is shown to be simple under general conditions (proof in Appendix A).    The final result of the paper is the following theorem where the automorphism group is thought of as a topological group with the permutation topology inherited from the action on the vertex set of the graph.

\medskip

{\bf Theorem~\ref{TPrimitiveSimple}}\ \ 
{\em Let $X$ be a locally finite connected primitive graph (meaning that the automorphism group is transitive and the automorphism group preserves no non-trivial proper equivalence relation on the vertex set) with infinitely many ends.  Then $G=\autx$ has an open topologically simple subgroup of finite index.} 

\medskip

A crucial part in the proof of this theorem is the fact that if $G$ is the automorphism group of a locally finite primitive graph with infinitely many ends then $G$ is non-discrete.  This follows from \cite[Theorem~2.5]{Smith2010} as is explained in Appendix B.

\section{General background and terminology}

\subsection{Graphs}

  We think of a {\em graph} $X$ as a pair $(VX,EX)$ where $VX$ is the vertex set 
and $EX$ is the set of edges.  
An edge is a two
element subset of $VX$.  If $e=\{x,y\}$ is an edge we say that $x$ and $y$ are {\em adjacent} and that $x$ and $y$ are the {\em end-vertices} of $e$.

The degree of a vertex  $x$ is the number of vertices adjacent to $x$.
A graph is said to be {\em locally finite\/} if all its vertices have 
finite degree.   A {\em path} in $X$ is a sequence
$v_0, v_1, \ldots, v_n$ of distinct vertices such that  $v_i$ is adjacent to 
$v_{i+1}$ for all $i=0, 1, \ldots, {n-1}$.  If the condition that the vertices be
distinct is dropped then we speak of a {\em walk}.   
A {\em ray\/} (also called a {\em half-line\/}) in a graph $X$ is a sequence 
$\{v_i\}_{i\in {\bf N}}$ of distinct vertices such that $v_i$ is adjacent to 
$v_{i+1}$ for all ${i\in \N}$.  A {\em line\/} (also called a {\em
  double ray\/})  in $X$ is a sequence 
$\{v_i\}_{i\in {\bf Z}}$ of distinct vertices such that $v_i$ is adjacent to 
$v_{i+1}$ for all ${i\in \Z}$.  For vertices $u$ and $v$ such that there is a path starting with $u$ and ending with $v$ we let $d(u,v)$ denote the minimum length
of such a path.  A path of length $d(u,v)$ starting with $u$ and ending with $v$ is called a {\em geodesic}.  
If $X$ is connected then $d$ is 
a metric on $VX$.  

The notion of a {\em graph} described above is often called a {\em simple graph} where an edge always has two distinct end-vertices and an edge is completely determined by its end-vertices.    When discussing quotients of graphs with a group action we need the more general concept of a {\em multigraph}.   A multigraph $X$ is a pair $(VX, EX)$ together with a map $t$ defined on the set $EX$ of edges such that the values of $t$ are either single vertices from $VX$ or two element subsets of $VX$ that represents the end vertices of $e$.
 
An {\em end} of a graph $X$ is defined as an equivalence class of rays such that two rays are said to be equivalent if there is a third ray that contains infinitely many vertices from both of them.  A connected graph $X$ has more than one end if and only if there is a finite set of vertices $F$ such that $X\setminus F$ has two distinct components that both contain rays.  If $X$ is a tree then two rays belong to the same end if and only if their intersection is a ray.  For more information about this concept consult \cite[Chapter 8]{Diestel2006} and \cite{Moller1995}.
 
\subsection{Permutation groups}

Let $G$ be a group acting on a set $Y$.  For $x\in Y$  let $G_x$ 
denote the {\em stabilizer of\/} $x$ in $G$; that is, $G_x$ is the subgroup of 
all elements in $G$ that fix $x$.  
For a subset $A$ of $Y$, define 
$$G_{(A)}=\{g\in G\mid g(a)=a\mbox{ for all }a\in A\}$$ 
and 
$$G_{\{A\}}=\{g\in G\mid gA=A\}.$$
The group $G_{(A)}$ is called the {\em pointwise stabilizer of} $A$ and the
  group $G_{\{A\}}$ the {\em setwise stabilizer of} $A$.
Two points $x, y$ are said to be in the same orbit of $G$ if there is an element $g\in G$ such that $g(x)=y$.  If any two elements in $Y$ are in the same orbit then we say that $G$ is transitive on $Y$.  The orbits of a stabilizer $G_x$ of a point $x\in Y$ are called the {\em suborbits} of $G$.
  
An action of a group $G$ on a set $Y$ defines a homomorphism from $G$ to the group $\sym Y$ of all permutations of $Y$.  If this homomorphisms is injective, i.e.\ the only element of $G$ that fixes all the points in $Y$ is the identity,  we say the action is {\em faithful}.  Then we can think of $G$ as a subgroup of $\sym Y$ and speak of $G$ as a {\em permutation group}.  

The automorphism group of a graph $X$ is denoted by \autx, and we think of \autx\ 
primarily as a 
permutation group on $VX$.  If \autx\ acts transitively on $VX$ then 
the graph $X$ is said to be {\em transitive}.  Varying slightly from the terminology above we define {\em the stabilizer of an edge} in a graph $X$ to be the subgroup fixing both end vertices of the edge.  

If a group $G$ acts on a set $Y$ we can construct a graph $X$ with vertex set $Y$ such that the action of $G$ on $Y$ gives an action on $X$ by automorphisms.  This is done by insisting that the edge set $EX$ is a union of orbits of $G$ on the set of two element subsets of $Y$.   Such graphs are called {\em (undirected) orbital graphs}.

A group is said to act
primitively on a set $Y$ if  the only $G$-invariant equivalence
relations on $Y$ are the trivial one (each equivalence class contains
only one element) and the universal one (there is only one equivalence class).
If a group $G$ acts transitively on a set $Y$ the following conditions are
equivalent:

{\bf (i)}  $G$ acts primitively.

{\bf (ii)}  if $y$ is a point in $Y$ then the subgroup $G_y$ is a
maximal subgroup of $G$.

{\bf (iii)} for every pair $x$ and $y$ of distinct points in $Y$ the orbital
graph with edge set $G\{x, y\}$ is connected.  

We say that a graph
$X$ is {\em primitive} if \autx\ acts primitively on $VX$.

\subsection{The permutation topology}

Let $G$ be a group acting on a set $Y$.  
The {\em permutation topology} on $G$ is defined
by choosing as a neighbourhood basis of the identity the family of
pointwise stabilizers of finite subsets of $Y$, i.e.\ a neighbourhood
basis of the identity is given by the family of subgroups
\[\{G_{(F)}\mid F\mbox{ is a finite subset of }Y\}.\]
For an introduction to the permutation topology see \cite{Moller2010}.

From this definition it is apparent that a sequence $(g_i)_{i\in \N}$ of elements in $G$ has an element $g\in G$ as a limit if and only if for every point
$y\in Y$ there is a number $N$ (possibly depending on $y$) such
that $g_n(y) =g(y)$ for every $n\geq N$.  
One could also use the property above describing convergence of sequences as a definition of the topology and think of the permutation topology as the topology of pointwise convergence.  If we think of $Y$ as having the discrete topology and elements of
$G$ as maps $Y\rightarrow Y$, then the permutation topology is equal to the compact-open topology.

A subgroup $U$ of $G$ is open if and only if there is a finite subset $F$ of $Y$ such that $G_{(F)}\subseteq U$.  One can also note that if $G$ is a permutation group on $Y$ then the permutation topology makes $G$ totally disconnected.  
Compactness has a natural
interpretation in the permutation topology.
A subset of a topological space is said to be {\em
 relatively compact} if it has compact closure.

\begin{lemm}{\rm (\cite[Lemma 1 and Lemma 2]{Woess1992}, cf. \cite[Lemma~2.2]{Moller2010})}\label{LCompact}
Let  $G$ be a group acting transitively on a set $Y$.
Assume that  $G$ is closed in the permutation topology and that all suborbits are finite.

(i) The stabilizer $G_y$ of a point $y \in Y$ is
compact.

(ii) A subset $A$ of $G$ is
 relatively compact
 in $G$ if and only if the set $Ay$ is finite for every point $y$
in $Y$.

Furthermore, if $A$ is a subset of $G$ and $Ay$
is finite for some $y\in Y$ then $Ay$ is finite for every
$y$ in $Y$.
\end{lemm}

A subgroup $H$ in a topological group $G$ is said to be cocompact if
$G/H$ is a compact space.  This concept has also a natural
interpretation in terms of the permutation topology.

\begin{lemm}\label{Lcocompact}
{\rm (\cite[Proposition~1]{Nebbia2000}, cf.~\cite[Lemma~2.3]{Moller2010})}
Let $G$ be a group acting transitively on a set $Y$. Assume
that $G$ is closed in the permutation topology and all suborbits are
finite.
Then a subgroup $H$ of $G$ is
cocompact if and only if $H$ has finitely many orbits on
$Y$.
\end{lemm}

\section{Groups acting on trees}

\subsection{Preliminaries on trees and group actions on trees}

In the present context {\em a tree} is a connected graph
that has no non-trivial cycles (i.e.\ there are no walks $v_0, v_1, \ldots,
v_n$ such that $v_0=v_n$ and $v_0, v_1, \ldots, v_{n-1}$  are distinct
and $n\geq 3$).    

In \cite{Tits1970} Tits classifies the automorphisms of a tree $T$.  First there are automorphisms that fix some vertex of $T$, then there are automorphisms that leave some edge of $T$ invariant but transpose its  end-vertices and finally there are translations.  An automorphism $g$ of $T$ is called a {\em translation} if there is a line $L=\{v_i\}_{i\in {\bf Z}}$ that is invariant under $g$ and there is a non-zero integer $k$ such that $g(v_i)=v_{i+k}$ for all $i\in\Z$.  The line $L$ is called the {\em axis} of the translation.  Suppose $e$ is an edge and $T_1$ is one of the components of $T\setminus e$.  If $g$ is an automorphism of $T$ such that $g(T_1)$ is a proper subset of $T_1$ then $g$ is a translation and the edge $e$ lies on the axis of $g$.  

When $L$ is a path (finite or infinite) in a tree
there is a well-defined map from the vertex set of the tree to the
vertex set of $L$ such that a vertex $x$ is mapped to the unique vertex in
$L$ that is closest to $x$.  This map will be denoted
with $\pr_L$.  For a vertex $x$ in $L$ the set $\pr_L^{-1}(x)$ is the vertex set of 
a subtree of $T$ which we call {\em the branch of} $T$ {\em at} $x$
(relative to $L$).  Let  $G$ be a group acting on $T$.  Note that the set $\pr_L^{-1}(x)$ is invariant under the group $G_{(L)}$.  Define
$G_{(L)}^x$ as the 
permutation group that we get by restricting the action of $G_{(L)}$
to $\pr_L^{-1}(x)$.  From the maps $G_{(L)}\rightarrow G_{(L)}^x$ we get a
homomorphism from the group $G_{(L)}$ to the
group $\prod_{x\in L} G^x_{(L)}$.  Following Tits in 
\cite{Tits1970} we say that a group $G$ acting on a tree $T$ has property
P if the homomorphism $G_{(L)}\rightarrow \prod_{x\in L} G_{(L)}^x$ above is an
  isomorphism for every path $L$ in $T$.  In \cite{CapraceDeMedts2011} property P is called {\em Tits' independence property}.  The essence of property P is that $G_{(L)}$ acts on each branch of $T$ at $L$ independently  of how it acts on the other branches.  In his groundbreaking paper
  Tits then goes on to prove \cite[Th\'eor\`eme~4.5]{Tits1970} that 
if a group $G$ acts on a tree $T$ such that property P is satisfied
  and no proper non-empty
  subtree is invariant under $G$ and no end of $T$ is fixed by $G$ 
then the subgroup $G^+$ of $G$
  generated by the stabilizers of edges is simple.  (Recall that the stabilizer of an edge $e=\{x,y\}$ in $G$ is defined as the subgroup fixing both $x$ and $y$.)

\subsection{Properties P, E and H}

A subtree $T'$ of a tree $T$ is called a {\em half-tree} of $T$ if $T'$ is one of the components of $T\setminus\{e\}$ for some edge $e$ in $T$ and both components of $T\setminus\{e\}$ are infinite.
For an edge $e=\{u, v\}$ in $T$ we let $T_{u, e}$ denote the component of $T\setminus \{e\}$ that contains $u$.  

A group acting on a tree $T$ is said to have property E if for every edge $e=\{u, v\}$ the stabilizer of $e$ acts independently on the two components of $T\setminus\{e\}$, i.e. $G_e=G_{(T_{u, e})}G_{(T_{v, e})}$.  If $G$ is a closed subgroup of $\aut T$ then properties P and E are equivalent, see \cite[Lemma 10]{Amann2003}. 

\begin{lemm}\label{LNoInvariant}
Let $G$ be a subgroup of $\aut T$ for some infinite tree $T$.

(i)  {\rm (\cite[Lemme 4.1]{Tits1970})}  The following two conditions are  equivalent:  (a) no proper non-empty subtree of $T$ is invariant under $G$  and (b)  for every vertex $x$ in $T$ the orbit $Gx$ intersects every half-tree of $T$.

Furthermore if (a) (and then (b) also) holds then the tree $T$ has no leaves (vertices of degree 1) and every edge defines two half-trees that both have unbounded diameter.

(ii)  {\rm (\cite[Lemme 4.4]{Tits1970})}  Suppose $N$ is a non-trivial subgroup of $\aut T$ normalized by $G$.   If no proper non-empty subtree of $T$ is invariant under $G$ and no end of $T$ is fixed by $G$ then the same is true about $N$.

(iii)  Suppose that no  proper non-empty subtree of $T$ is invariant under $G$.  If $e$ is an edge in $T$ then $G$ contains a translation $g$ such that $e$ is  in the axis of $g$.
\end{lemm}

{\em Proof.}  We only need to prove part (iii); parts (i) and (ii) are proved in \cite{Tits1970} except the addendum in part (i) which is obvious.   Let $u$ and $v$ denote the end vertices of $e$.  By part (i) we can find an element $g_1\in G$ such that $g_1(u)\in T_{v,e}$.  If $g_1(T_{v,e})\subsetneq T_{v,e}$ then $g_1$ is the translation we are seeking. If not we find an element $g_2\in G$ such that $g_2(v)\in  T_{u,e}$.  As before, if 
$g_2(T_{u,e})\subsetneq T_{u,e}$ then $g_2$ is the translation we are seeking but if not then $g=g_1g_2^{-1}$ is the translation we are after.  The result in part (iii) and the argument in the proof are well-known and can for example be seen in a more general context in \cite{Jung1981}.

\bigskip

\begin{lemm}\label{LPropertyH}
Let $T$ be a tree and $G$ a subgroup of $\aut T$.  Assume that no proper non-empty subtree of $T$ is invariant under $G$.   

(i)  Suppose that there is some edge $e=\{u,v\}$ in $T$ such that the pointwise stabilizers of both the half-trees $T_{u,e}$ and $T_{v,e}$ are trivial.  Then the stabilizer of every half-tree in $T$ is trivial.

(ii)  Suppose that there is some edge $e=\{u,v\}$ in $T$ such that the pointwise stabilizers of both the half-trees $T_{u,e}$ and $T_{v,e}$ are non-trivial.  Then the stabilizer of every half-tree in $T$ is non-trivial.

(iii)  Suppose that there is some edge $e=\{u,v\}$ in $T$ such that the pointwise stabilizer of $T_{u,e}$ is trivial but the pointwise stabilizer of  $T_{v,e}$ is non-trivial.  Then $G$ must fix an end of $T$. 
\end{lemm}

{\em Proof.}   (i) Let $f$ be an edge in $T$.  This edge defines two half-trees of $T$ which we denote with $T_1$ and $T_2$.  By Lemma~\ref{LNoInvariant}(i) there is an element $g\in G$ such that $g(e)\in T_1$.  Then either $g(T_{u,e})\subseteq T_1$ or $g(T_{v,e})\subseteq T_1$.  If we assume, for instance, that $g(T_{u,e})\subseteq T_1$ then 
$G_{(T_1)}\subseteq G_{(g(T_{u,e}))}=gG_{(T_{u,e})}g^{-1}=\{1\}$.  In the same way we show that $G_{(T_2)}$ is trivial.

(ii)  Let $f$, $T_1$ and $T_2$ be as above.  Find an element $g\in G$ such 
that $g(e)\in T_1$.  Then either  $g(T_{u,e})\subseteq T_1$ or  $g(T_{v,e})\subseteq T_1$.  Say, for the sake of the argument, $g(T_{u,e})\subseteq T_1$ and then $T_2\subseteq g(T_{v,e})$.  Then $f\in T_{v,e}$ and $T_2\subseteq g(T_{v,e})$.  Hence 
$G_{(T_2)}\supseteq G_{(gT_{v ,e})}=gG_{(T_{v ,e})}g^{-1}\neq\{1\}$.  So $G_{(T_2)}$ is non-trivial. In the same way we show that  $G_{(T_1)}$ is also non-trivial.

(iii)  Looking at parts (i) and (ii) we conclude that it is true for every edge $f=\{w,z\}$ in $T$ that the pointwise stabilizer of one of the half-trees defined by $f$ is trivial and the pointwise stabilizer of the other one is non-trivial.   If the pointwise stabilizer of $T_{w,f}$ is trivial then we think of $f$ as an directed arc with initial vertex $w$ and terminal vertex $z$.  The edge $e=\{u,v\}$ is oriented so that $u$ is the initial vertex and $v$ the terminal vertex.  Do this for every edge in $T$ and note that the direction of the edges is preserved by the action of $G$.    

Consider an edge $f=\{w,z\}$ that is contained in the half-tree $T_{v,e}$.  Assume that $w$ is closer to $v$ than $z$.  Then $T_{w,f}\supseteq T_{u,e}$ and thus $G_{(T_{w,f})}\subseteq G_{(T_{u,e})}=\{1\}$.  Hence the initial vertex of $f$ is $w$ and the terminal vertex is $z$.  Another way to describe this is to say that the edge $f$ is directed away from $v$.  We see that every vertex in $T_{v,e}$ is the terminal vertex of precisely one directed edge.  The half-tree $T_{v,e}$ contains vertices from every $G$-orbit on the  vertex set of $T$ and thus it is true for every vertex in $T$ that it is the terminal vertex of precisely one directed edge.  Let $R_1=v_0, v_1, \ldots$ be a ray in $T$ such that each edge $\{v_i, v_{i+1}\}$ is directed so that $v_i$ is the terminal vertex.  Suppose now that $R_2=w_0, w_1, w_2, \ldots$ is another such ray in $T$.  If the two rays intersect then they both belong to the same end.  Suppose the  two rays are disjoint.  Select a path  $u_0, u_1, \ldots, u_n$ of shortest possible length such that $u_0$ is a vertex in $R_1$ and $u_n$ is a vertex in $R_2$.  Then $u_1$ is not in $R_1$ and thus the edge $\{u_0,u_1\}$ is directed so that $u_0$ is the initial vertex.  Note also that $u_1$ is not on $R_1$. Similarly the edge $\{u_{n-1}, u_n\}$ has $u_n$ as an initial vertex.  Now we see that there must be a number $k$ such that the edges $\{u_{k-1}, u_k\}$ and $\{u_{k}, u_{k+1}\}$ both have $u_k$ as a terminal vertex.  This is a contradiction and thus the two rays $R_1$ and $R_2$ belong to the same end $\omega$, which is clearly fixed by $G$.

\begin{defi}
A group $G$ that is a subgroup of the automorphism group of some tree $T$ is said to have property H if the pointwise stabilizer of every  half-tree in $T$ is non-trivial. 
\end{defi}

Let $G^{++}$ denote the closure of the subgroup generated by the pointwise stabilizers of all the half-trees in T.  Clearly $G^{++}$ is normal in $G$.
If $G$ has no proper non-empty invariant subtree and  does not fix an end of $T$ then  Lemma~\ref{LPropertyH} above shows that if $G$ does not have property H then the pointwise stabilizer of every half-tree is trivial and then
property H is equivalent to the property that the group $G^{++}$ is non-trivial.
Note also that if $G$ has property H then $G$ is not discrete.

The relationship between properties P and H is not simple.  In the case that $G$ is not discrete and has no non-empty proper invariant subtree then property P implies property H.  On the other hand the example below shows that property H does not imply property P.

\medskip

{\em Example.}  Let $T$ be a tree and $f:VT\rightarrow I$ some map defined on the vertex set of $T$.  In his paper Tits \cite{Tits1970} studies the group  $\aut\!\!_f\; T=\{g\in \aut T\mid f\circ g=f\}$.  This group clearly has property P.  One can also study the group $G$ of all automorphisms of $T$ that preserve the equivalence relation defined by the fibers of $f$.  It is not to be expected that this group has property P, but in many cases it will have property H.  

Consider the case of a regular tree $T$ of degree 6.  Colour all the vertices in one part of the natural bipartition red and then colour the vertices in the other part of the natural bipartition with three different colours so that each red vertex is adjacent to two vertices of each colour.  The group $C$ of automorphisms of $T$ that map every vertex to a vertex of the same colour has property P and is simple by Tits' theorem \cite[Th\'eor\'eme~4.5]{Tits1970}.  The group $G$ of automorphisms that leave the partitioning of the vertices given by this colouring invariant does not have property P but it has property H and it is easy to see that $G^{++}=C$.

\begin{theo}\label{TSimple}
Let $G$ be a closed subgroup of the automorphism group of some tree $T$.   Assume that no proper non-empty subtree of $T$ is invariant under $G$ and no end of $T$ is fixed by $G$.  Assume also that $G$ has property H.
If $N$ is a non-trivial closed subgroup of $G$ normalized by $G^{++}$ then $N$ contains $G^{++}$.   In particular, the subgroup $G^{++}$ is topologically simple.  
\end{theo}

{\em Proof.}  First note that by Lemma~\ref{LNoInvariant}(ii) we see that $G^{++}$ does not leave any proper non-empty subtree invariant and $G^{++}$ does not fix an end.   Now we apply Lemma~\ref{LNoInvariant}(ii) again but this time to $G^{++}$ and $N$ and find out that $N$ does not leave any proper non-empty subtree invariant and does not fix an end.
Let $e=\{u,v\}$ be an edge in $T$.  
By part (iii) of Lemma~\ref{LNoInvariant} we see that there is a translation $h\in N$ such that $h(T_{v,e})\subsetneq T_{v,e}$ and $T_{u,e}\subsetneq h(T_{u,e})$.  Suppose $g\in G_{(T_{u,e})}$.  Set $f_n=gh^ng^{-1}h^{-n}$.  
Since $N$ is  normalized by $G^{++}$ we see that $f_n\in N$ for every $n\geq 0$.  The element $h^ng^{-1}h^{-n}$ fixes the half-tree $h^{n}(T_{u,e})$ and in particular $T_{u,e}$ is fixed by $f_n$.  If we consider $T_{v,e}\setminus h^{n}(T_{v,e})$ then this part of the tree is fixed by $h^ng^{-1}h^{-n}$ and thus $f_n$ acts on this part like $g$.  Whence we see that $f_n\rightarrow g$ when $n\rightarrow \infty$.  From this argument we conclude that $G_{(T_{u, e})}$ is contained in $N$.  Of course one can apply the same argument to show that $G_{(T_{v, e})}$ is contained in $N$.  We conclude that $G^{++}$ is contained in $N$.  Now it is clear that $G^{++}$ is topologically simple.

\bigskip

{\em Remark.}  The contraction group for an automorphism $\alpha$ of a topological group $G$ is defined as the subgroup of all elements $g\in G$ such that $\alpha^n(g)\rightarrow 1$, see \cite{BaumgartnerWillis2004}.   In the above proof $h^ng^{-1}h^{-n}\rightarrow 1$ and hence $g^{-1}$ belongs to the  contraction group for the inner automorphism of $G$ defined by $h$.

\begin{coro}  
Let $G$ be a closed subgroup of $\aut T$ for some tree $T$.  Suppose that $G$ has property H and does not stabilize a proper non-empty subtree or fix an end.   Then $G^{++}$ is the unique minimal closed normal subgroup of $G$.
\end{coro}  

The {\em quasi-center} $QZ(G)$  of a topological group $G$ consists of all elements with an open centralizer.  Caprace and De Medts show in \cite[Proposition~3.6]{CapraceDeMedts2011} that a closed subgroup $G$ of the automorphism group of some tree $T$ such that $G$ satisfies property P and does not have any proper non-empty invariant subtrees has a trivial quasi-center.   A simple adaptation of their proof gives an analogous result for groups with property H.

\begin{prop}
Let $G$ be a closed subgroup of $\aut T$ for a tree $T$.  Assume that $G$ has property H and that $G$ leaves no proper non-empty subtree of $T$ invariant.  Then the quasi-center of $G$ is trivial.  
\end{prop}
 
{\em Proof.}  Suppose $g$ is an element of $G$ that has an open centralizer.  Let $v$ be a vertex of $T$.  Then there is a finite set $S$ of vertices such that $G_{(S)}$ is contained in the centralizer of $g$.  If necessary we can replace $S$ with $S\cup\{v\}$ so we can assume that $v\in S$ and we may also safely assume that $S$ is a subtree of $T$.  Let $\tilde{S}$ be the subtree of $T$ containing every vertex of $T$ that is fixed by $G_{(S)}$.  
Because $g$ centralizes $G_{(S)}$ the tree $\tilde{S}$ is invariant under $g$.    Suppose $e=\{u,w\}$ is an edge in $T$ such that the vertex $u$ is in  $\tilde{S}$ but $w$ is not.  Using property H we can find a nontrivial element $h\in G_{(T_{u,e})}$.    But $\tilde{S}\subseteq  T_{u,e}$ so $h\in G_{(S)}$ and $g$ commutes with $h$ and $ghg^{-1}=h$.  Note that $h=ghg^{-1}$ fixes $g(T_{u,e})=T_{g(u),g(e)}$.  Since $g(\tilde{S})=\tilde{S}$, we see that if $g(u)\neq u$ then $T_{w,e}\subseteq g(T_{u,e})$ and then $h$ would fix pointwise both $T_{u,e}$ and $T_{w,e}$ -- a contradiction.    Hence we conclude that $g(u)=u$.   Therefore $g$ fixes every vertex in $\tilde{S}$ that is adjacent to some vertex not in $\tilde{S}$.  Suppose now that $e=\{u,w\}$ is an edge in $T$ such that $u$ is in $S$ but $w$ is not in $S$.  If the edge $e$ is not in $\tilde{S}$ then the above argument shows that $g$ fixes $u$.  On the other hand, if $e$ is in $\tilde{S}$ then $G_{(T_{u,e})}\neq\{1\}$ and $G_{(T_{u,e})}\subseteq G_{(S)}$ and thus $G_{(S)}$ moves some vertex in $T_{w,e}$.  Therefore $T_{w,e}$ is not contained in $\tilde{S}$.  From this we infer that $T_{w,e}$ contains a vertex $z$ in $\tilde{S}$ that is adjacent to a vertex not in $\tilde{S}$ and thus $z$ is fixed by $g$.  Applying this argument to every edge with precisely one of its end vertices in $S$ and we conclude that $g$ must fix a vertex in every component of $T\setminus S$.  Since $S$ is finite we now see that every vertex in $S$ is fixed by $g$ and in particular $g$ fixes $v$.  Since $v$ was arbitrary we conclude that that $g$ fixes every vertex of $T$ and that $g=1$.

\bigskip

The following is an analogue of Proposition~3.8 from the paper \cite{CapraceDeMedts2011} of Caprace and De Medts and the proof uses the same argument.     

\begin{prop}\label{PSubtree}{\rm (Cf.~\cite[Proposition~3.8]{CapraceDeMedts2011})}
Let $G$ be a closed subgroup of the automorphism group of some tree $T$ that leaves no proper non-empty subtree invariant.  
Suppose $H$ is a non-compact open subgroup of $G$ that does not fix an end of $T$ and $T'$ is a minimal invariant subtree for $H$.   Then for every edge $e$ in $T'$ the group $H$ contains the pointwise stabilizers in $G$ of the two half-trees of $T$ defined by the edge $e$. 
\end{prop}

{\em Proof.}  Since $H$ is non-compact the tree $T'$ is infinite and every edge $e$ in $T'$ splits $T'$ up into two half-trees.
From Lemma~\ref{LNoInvariant}(iii) above we see that for every edge $e$ in $T'$ there is a hyperbolic element $h$ in $H$ such that $e$ is on the axis of $h$.  Since $H$ is open there is a finite set of vertices such that the pointwise stabilizer $G_{(S)}$ is contained in $H$.  Let $T_1$ and $T_2$ denote the two half-trees of $T$ defined by $e$.  Using a suitable power of $h$ we can assume that $h^n(S)\subseteq T_1$ and then $G_{(T_1)}\subseteq h^nG_{(S)}h^{-n}\subset H$.   Similarly we can show that $G_{(T_2)}\subseteq H$.  

\bigskip

If $H$ is compact then the tree $T'$ either has just a single vertex or consists of an edge with its end vertices and $H$ then contains an element that transposes the two end vertices.  If it is assumed that $G$ is topologically simple (like in \cite[Proposition~3.8]{CapraceDeMedts2011}) then 
$G$ acts without inversion on $T$ and the latter possibility above can not occur and the conclusion of Proposition~\ref{PSubtree} holds trivially.

\begin{lemm}{\rm (Cf.~\cite[Lemma~3.11]{CapraceDeMedts2011})}\label{LNoFree}
Let $G$ be a closed edge transitive subgroup of the automorphism group of some tree $T$.  If $G$ is simple and has property H then there is no vertex $v$ in $T$ such that the action of $G_v$ on the set of edges with $v$ as an end vertex is free.
\end{lemm}

The argument in the proof of lemma above is the same as in \cite{CapraceDeMedts2011}.

The following is an version of \cite[Theorem A]{CapraceDeMedts2011}, but here  
property H is assumed instead of property P.  Caprace and de Medts derive their theorem from \cite[Theorem ~3.9]{CapraceDeMedts2011} and their argument also works for this version where we use Proposition~\ref{PSubtree} instead of \cite[Proposition~3.8]{CapraceDeMedts2011} and Lemma~\ref{LNoFree} instead of 
\cite[Lemma~3.11]{CapraceDeMedts2011}.

\begin{theo}\label{TTheorem A}{\rm (Cf.~\cite[Theorem A and Theorem 3.9]{CapraceDeMedts2011})}
Let $T$ be a tree all of whose vertices have degree at least 3.  Suppose $G$ is a topologically simple closed subgroup of $\aut T$ which does not stabilize any proper non-empty subtree and which satisfies property H.  Then the following conditions are equivalent.

(i) Every proper open subgroup of $G$ is compact.

(ii)  For every vertex $v\in VT$, the induced action of $G_v$ on the edges that have $v$ as an end vertex is primitive.  In particular the action of $G$ on the set of edges of $T$ is transitive. 
\end{theo}

\section{The automorphism group of a graph with\\ connectivity 1}\label{SConnectivity}

A connected graph $X$ is said to have {\em connectivity} 1 if there is a vertex $x$ in $X$ such that $X\setminus x$ is not connected.  Such a vertex $x$ is called a {\em cutvertex}.  If a transitive
graph has a cutvertex then every vertex is a cutvertex. 

The {\em blocks} (called {\em lobes} in \cite{JungWatkins1977}) 
of a graph $X$ with connectivity 1 are
the maximal connected subgraphs that do not have connectivity 1.  

In this section, we obtain some simplicity results on the automorphism group of a transitive graph $X$ with connectivity 1. We use Tits' simplicity theorem \cite[Th\'eor\`eme~4.5]{Tits1970}.

From a graph $X$ with connectivity 1 we can construct a tree $T_X$
called the {\em block graph} of $X$.  The vertex set of $T_X$
is the union of the set of blocks of $X$ and the set of cutvertices in
$X$.  The set of edges in $T$ consists of all pairs $\{x, B\}$ where $x$ is a
cutvertex and $B$ a block and $x$ is in $B$.  
The set of cutvertices and the
set of blocks thus form the parts of the natural bipartition of the tree $T_X$.
The automorphism group of $X$ acts on $T_X$.  

\begin{lemm}\label{LPropertyP}
Let $X$ be a transitive graph with connectivity 1.  
The action of $G=\autx$ on $T_X$ has property P.
\end{lemm}

{\em Proof.}  This can be seen directly or by noting that the action
of $G$ on $T_X$ clearly has property E, and as $G$ is a closed
permutation group, the action has property P.  

\bigskip

This lemma allows us to prove certain simplicity results for $\autx$. We say that a group $G$ acting on a set $Y$ is generated by stabilizers of points if the stabilizers in $G$ of the points in $Y$ generate $G$.

\begin{theo}\label{TConnectivity1Simple}
Let $X$ be a transitive graph with connectivity $1$ and $G=\autx$. Let $n$ 
be the number of blocks a vertex in $X$ lies in.

(i) If the automorphism group of some block is not transitive, then $G$ is not simple.

(ii) If the automorphism group of every block is transitive and generated by vertex stabilizers, then $G$ is simple, unless $n=2$ and any two blocks are isomorphic, in which case $G$ has a normal simple subgroup of index $2$.
\end{theo}

\begin{coro}\label{CPrimSimple}
Let $X$ be a primitive graph with connectivity $1$.  If each vertex is
contained in more than two blocks then the group
$G=\autx$ is simple.   If each vertex is only contained in two blocks
then $G$ has a simple normal subgroup of index $2$.
\end{coro}

The proofs of Theorem~\ref{TConnectivity1Simple} and Corollary~\ref{CPrimSimple} can be found in Appendix A together with the necessary background.

\section{Automorphism groups of primitive graphs with infinitely many ends}

\begin{theo}\label{TPrimitiveSimple}
Let $X$ be a locally finite connected primitive graph with infinitely many ends.  Then $G=\autx$ has an open topologically simple subgroup of finite index. 
\end{theo}

{\em Proof.}  {\bf Step 1} in the proof is to define an action of $G$ on a tree $T$.    
It is shown in  \cite[Proposition~3]{Moller1994} that if $X$ is a locally finite primitive graph with more than one end and $G$ is a group acting primitively on $X$ by automorphisms then there is a pair of vertices $x, y$ in $X$ such that the graph $Y$ with the same vertex set as $X$ and edge set $EY=G\{x, y\}$ is connected and has connectivity 1.  Note that the action of $G$ on the vertex set of $X$ (the same as the vertex set of $Y$) gives an action of $G$ by automorphisms on $Y$.  The group $G$ now acts on the block graph $T_Y$ that is a tree, as explained in Section~\ref{SConnectivity}.

{\bf Step 2 }is to show that the action of $G$ on $T_Y$ is faithful, fixes no end of $T_Y$ and leaves no proper non-empty subtree invariant.    As explained in Section~\ref{SConnectivity} we can think of a vertex in $X$ also as a vertex in $T_Y$ and identify the vertex set of $X$ with one of the parts of the natural bipartition of the vertex set of $T_Y$.  The action of $G$ on $T_Y$ is thus obviously faithful and there is no proper non-empty invariant subtree.  Suppose that $G$ fixes some end of $T_Y$.  We want to define a $G$ invariant proper non-trivial equivalence relation on the vertex set of $X$ contradicting the assumption that $G$ acts primitively on $X$. 
Take a ray  $R=v_0, v_1, v_2, \ldots$ in $T$ belonging to an end $\omega$ fixed by $G$ and say that vertices $u$ and $v$ are related if there is a number $N(u,v)$ such that $d(u, v_i)=d(v,v_i)$ for all numbers $i$ larger than $N(u,v)$.  It is left to the reader to show that this is an equivalence relation and does not depend on the choice of the ray $R$.  The equivalence classes are often called {\em horocycles}.  Since the end $\omega$ is fixed by $G$ this equivalence relation is invariant under $G$.  Restricting to the vertex set of $X$ (which we think of as a subset of the vertex set of $T_Y$) we see that this would give a proper non-trivial $G$ invariant equivalence relation on the vertex set of $X$ contradicting the assumption that $G$ acts primitively on $X$.  Hence it is impossible that $G$ fixes an end of $T_Y$.  

{\bf Step 3 }is to show that the action of $G$ on $T_Y$ has property H.    An edge $\{x, B\}$ in $T_Y$ where $x$ is a vertex in $Y$ and $B$ is a block in $Y$  gives a partition of $T_Y$ into two half-trees and that in turns gives a partition of the vertex set of $Y$ into two disjoint parts $C_{x}$ and $C_{B}$.  Let $S_Y$ be the set of all the edges in the block $B$  that have $x$ as a endvertex.  If we remove the edges in $S_Y$ from $Y$ then we get a graph with two components that have vertex sets $C_{x}$ and $C_{B}$, respectively.  

Define now $S_X$ as the set of edges in $X$ that have one endvertex in $C_{x}$ and the other one in $C_{B}$.  Because $X$ is locally finite and $G$ is transitive on the vertex set of $X$ we see that $G$ has only finitely many orbits on pairs $\{u, v\}$ of adjacent vertices in $X$.  The action of $G$ on the vertex set of $X$ (the same as the vertex set of $Y$) induces automorphisms of both $X$ and $Y$ and we see that there is a constant $k$ such that if $u$ and $v$ are adjacent vertices in $X$ then $d_Y(u,v)\leq k$.  
The graph $Y$ is locally finite and thus there are only finitely many pairs of vertices $u\in C_{x}$ and $v\in C_{B}$ such that $d_Y(u,v)\leq k$.  Now it follows from the above that the set $S_X$ is finite.

Define $H$ as the subgroup of $G$ consisting of all the elements of $G$ that fix all the edges in $S_X$ and their endvertices.  Since the set $S_X$ is finite, the group $H$ is open in the permutation topology on $G$.  This groups leaves the sets $C_{x}$ and $C_B$ invariant.  
It follows from \cite[Theorem~2.5]{Smith2010} that the group $G$ in the permutation topology is non-discrete (see Appendix B for a detailed explanation).  
The set $S_X$ separates the sets $C_{x}$ and $C_B$ (i.e. any path between an vertex in $C_{x}$ and a vertex in $C_B$ contains an edge from $S_X$).  Because $G$ is the full automorphism group of $X$ then $H$ acts independently on $C_{x}$ and $C_B$, i.e. $H=H_{(C_x)}H_{(C_B)}$.  As $H$ is nontrivial, $H_{(C_x)}$ or $H_{(C_B)}$ is non-trivial.  But $H_{(C_x)}$ is contained in $G_{(T_{x,e})}$ and  $H_{(C_B)}$ is contained in $G_{(T_{B,e})}$.  Hence the action of $G$ on $T_Y$ has property H.  

The final {\bf step} is an application of Theorem~\ref{TSimple}.
As stated above the group $G^{++}$ is in this case a non-trivial topologically simple open normal subgroup of $G$ and every closed non-trivial normal subgroup of $G$ contains $G^{++}$.    Since $G$ acts primitively on the vertex set of $X$ the normal subgroup $G^{++}$ acts transitively on the vertex set of $X$.
By Lemma~\ref{Lcocompact} we conclude that $G^{++}$ is cocompact in $G$.  Because $G^{++}$ is open we know that the quotient space $G/G^{++}$ is discrete and since it is also compact we see that it must be finite.  Hence $G^{++}$ has finite index in $G$. 

\bigskip

{\em Remark.}  The result \cite[Proposition~3]{Moller1994} about primitive graphs referred to above is proved by using the theory of structure trees developed and described for instance in \cite{DicksDunwoody1989}, \cite{Moller1995} and \cite{ThomassenWoess1993}.   Using this theory it is possible to apply Theorem~\ref{TSimple} more generally to automorphism groups of locally finite graphs with infinitely many ends.

\section*{Appendix A:  Automorphism groups of graphs with connectivity 1}\label{AConnectivity1}

This appendix contains the proofs of Theorem~\ref{TConnectivity1Simple} and Corollary~\ref{CPrimSimple} together with necessary background discussion.

For a graph $X$ with connectivity 1 we let $B_i$, $i\in I$, denote
a family of representatives for the 
isomorphism types of blocks in $X$.   Furthermore, use $B_i^{(j)}$
for $j\in J_i$ to denote the orbits of $\aut B_i$ on the vertex set
of $B_i$.  For a vertex $x$ in $X$ we let $m_i^{(j)}(x)$ be the number of blocks of type $B_i$ that contain $x$ in the orbit $B_i^{(j)}$.  
Jung and Watkins in \cite[Theorem 3.2]{JungWatkins1977} show that a graph $X$ of connectivity 1 is transitive if and only if all the functions $m_i^{(j)}(x)$ are constant on $VX$.

Consider an action of a group $G$ on a tree $T$ and
assume that the action has property P.  Let $G^+$ denote the subgroup
of $G$ generated by the stabilizers of edges.  The subgroup $G^+$ is
simple by Tits' theorem \cite[Th\'eor\`eme~4.5]{Tits1970}.  To decide if $G$ is simple we must investigate when $G=G^+$.  
For a vertex $x$ in $T$ we define $T_1(x)$ as the set of vertices adjacent to $x$.  The set  $T_1(x)$ is invariant under $G_x$ and we define 
$G_x^{T_1(x)}$ to be the permutation group that $G_x$ induces on $T_1(x)$.

\begin{lemm}\label{LStabilisers}
If $x$ is a vertex in $T$ then the group $G^+$ contains the group
$G_{(T_1(x))}$.  If the group $F=G_x^{T_1(x)}$ is generated by
stabilizers of vertices (thought of as a permutation group on
$T_1(x)$) then $G^+$ contains $G_x$.  
\end{lemm} 

{\em Proof.} The first part of the Lemma is obvious because if $e$ is
an edge in $T$ with $x$ as an end-vertex then $G_{(T_1(x))}$ is contained
in $G_e$ and thus $G_{(T_1(x))}$ is contained in $G^+$.

Consider now the action of $G_x^+$ on $T_1(x)$.  For a vertex $y$ in
$T_1(x)$ the group $G^+$ contains the stabilizer of the edge between
$x$ and $y$ and thus the stabilizer of $y$ in $G_x^{T_1(x)}$ is
contained in the permutation group induced by $G^+$ on $T_1(x)$.
Hence, if $F$ is generated by stabilizers then $G^+$ induces the full
group $F$ on $T_1(x)$.  Since $G^+$ contains the group $G_{(T_1(x))}$ we  conclude that $G^+$ contains the full stabilizer $G_x$.  

\begin{lemm}\label{LNotSimple}
Let $Y$ denote the quotient graph of $T$ by the action of $G$.  If the
graph $Y$ is not a tree then $G$ is not simple.   If $Y$ is a tree
then $G$ is generated by stabilizers of vertices.
\end{lemm} 

{\em Proof.} Let $R$ denote the normal subgroup of $G$ generated by all the stabilizes of vertices.  By \cite[Corollary 1 in \S
  5.4]{Serre1980} the quotient group $G/R$ is
isomorphic to the fundamental group of $Y$, and this is non-trivial if
$Y$ is not a tree.

\bigskip

{\em Proof of Theorem~\ref{TConnectivity1Simple}.}  (i) Suppose that $B$ is some block of $X$ such that the
automorphism group of $B$ is not transitive.  Say $x$ and $y$ are
vertices in $B$ that belong to different orbits of
$\aut B$.  The edges $\{x,B\}$ and  $\{y,B\}$ in
$T_X$ belong then to different orbits of $\autx$ and have therefore
different images in the quotient graph of $T_X$ under the action of
$G$.  But their images in the quotient graph have the same end-vertices and thus the quotient graph is not a tree.  By Lemma~\ref{LNotSimple}, $\autx$ is not simple.

(ii) The stabilizer in $G$ of a vertex $x$ in $X$ acts on the set of blocks that contain $x$ as a direct product of symmetric groups. If $n>2$ the permutation group induced by $G_x$ on the neighbours of $x$ in $T$ is generated by stabilizers, and therefore  $G_x \leq G^+$.

If $n=2$ and the two blocks containing a given vertex are not isomorphic, the same conclusion obviously holds since $G_x$ acts trivially on $T_1(x)$.

Let $B$ be a block of $X$ and think of $B$ as a vertex in $T_X$. The vertices in $T_X$ contained in $T_1(B)$ correspond to the vertices
in the block $B$ and $G_B$ induces the full automorphism group of $B$
on $T_1(B)$. Since $\aut B$ is generated by stabilizers, we have $G_B \leq G^+$.

The quotient graph $Y$ of $T_X$ by the action of $G$ has one vertex $\tilde{x} $ for the orbit of $G$ on
the vertices in $T_X$ corresponding to vertices in $X$ and one vertex
for each orbit on the blocks, joined to $\tilde{x} $.  This is a tree and by
Lemma~\ref{LNotSimple} we see that $G$ is generated by the stabilizers
of vertices and hence $G^+=G$.

\par Finally assume $n=2$ and all blocks are isomorphic. In this case each vertex in the tree $T_X$ corresponding to a vertex in $X$ has degree 2. Construct a new graph $T'_X$ such that the set of vertices is the set of blocks of $X$ and two vertices in $T'_X$ are adjacent if and only if the corresponding blocks have a common vertex. The condition that each vertex is contained in just two blocks guarantees that $T'_X$ is a tree. The assumption that all the blocks are isomorphic says that $G$ acts transitively on the vertex set of $T'_X$. One now sees that $G$ cannot be simple because
the subgroup $N$ of $G$ preserving the classes of the natural bipartition of $T'_X$ is normal in $G$ with index 2. It is clear that $N$ is generated by the stabilizers of vertices in $T'_X$ (i.e. the stabilizers in $G$ of the blocks of $X$). The argument above shows that the group $G^+$ contains all the stabilizers of blocks in $X$ and
thus $N=G^+$ and $N$ is simple.  
\bigskip

{\em Comment. } Assume $n=2$ and the two blocks containing a given vertex are not isomorphic. By the classification of Jung and Watkins of transitive graphs with connectivity 1 described above we see that the automorphism group of each block must be transitive. 

\bigskip

{\em Proof of Corollary~\ref{CPrimSimple}.}  
Jung and Watkins in \cite{JungWatkins1977} also give a complete description of primitive graphs with connectivity 1.  In these graphs each block is a
primitive graph, any two blocks are isomorphic and each block has at least three vertices. 

  Each block is a primitive graph (therefore, transitive) with at least three vertices and therefore the automorphism group of a block is generated by stabilizers.  (It is impossible that the automorphism group of a primitive graph is regular.)  Because all the blocks are isomorphic the result now follows from Theorem~\ref{TConnectivity1Simple}.

\section*{Appendix B:  Primitive graphs with infinitely many ends}\label{APrimitive}

The purpose of this appendix is to explain how Theorem~2.5 in Smith's paper \cite{Smith2010} implies that if a group acts primitively on an infinite locally finite graph with more than one end then the stabilizer of a vertex is infinite.  

We start by proving a results for group actions on trees and then use the tree described in Step 1 of the proof of Theorem~\ref{TPrimitiveSimple}
 to prove our results for a group acting primitively on a locally finite graph with infinitely many ends.

\begin{lemm}\label{LSimon}
Let $G$ be a group acting on a tree $T$.  Let $x$ and $y$ be distinct vertices
in $T$.  Assume that on the path between $x$ and $y$ is a vertex $z$, distinct from both $x$ and $y$, such that $G_{x, z}=G_{y, z}$.  Suppose that $d(x,z)\leq d(z,y)$.  If $h\in H=\langle G_x, G_y\rangle$ 
then either $h$ fixes $y$ or $d_T(y, h(y))>d(x,y)$.
\end{lemm}

{\em Proof.}  Set $A=G_x$,  $B=G_y$ and
$C=A\cap B=G_{x,y}$.  Let $\{a_i\}_{i\in I}$ be a set
  of coset representatives for $C$ in $A$ and, similarly, let
  $\{b_j\}_{j\in J}$ be a set of coset representatives for $C$ in $B$.
Assume that the identity element is included in
  both families.  

Set $k=d_T(x,y)$.  For $g\in A$ we define $x(g)$ as a
  vertex in $[x, y]$ that is fixed by $g$ and is in the greatest distance from $x$.  If $g\in B$  define
  $y(g)$ similarly.  For an element $g\in
  A\setminus C$ the condition in the lemma means that $d_T(x,x(g))<d_T(x,z)\leq k/2$
  and then $d_T(y,x(g))>d_T(y,z)\geq k/2$.  Similarly, if $g\in B\setminus C$ then $d_T(y,y(g))<d_T(y,z)$ and $d_T(x,y(g))>d_T(x,z)$.  
Recall that if $v$ is a vertex  in $T$ then $\pr_{[x,y]}(v)$ is defined as the vertex on the geodesic $[x,y]$
that is closest to $v$.

Write $h\in H$ as
  $h=b_{j_l}a_{i_l}\cdots b_{j_2}a_{i_2}b_{j_1}a_{i_1}c$ where $c\in C$ and none of the
  $a_i$'s and $b_j$'s is the identity element with the possible
  exceptions of $a_{i_1}$ and $b_{j_l}$.  We use induction over $l$.
Our induction hypothesis is that if $h$ is as above and $h$ does not fix $y$ then $d_T(y, h(y))>k$ and if 
$b_{j_l}\neq 1$ then $\pr_{[x,y]}(h(y))=y(b_{j_l})$.

To start with it is obvious that if $l=0$ or  $l=1$ and $a_{i_1}=1$ then $h$ fixes $y$.  Assume now
that $a_{i_1}\neq 1$.  Then $a_{i_1}$ is in $A\setminus C$ and does not fix $y$.  Recall that $d_T(y, x(a_{i_1}))>d_T(y,z)\geq k/2$.
The geodesic from $y$ to $a_{i_1}(y)$ goes
through the vertex $x(a_{i_1})$.  Hence $d_T(y, a_{i_1}(y))=2d_T(y,
x(a_{i_1}))>k$.   
 Because $b_{i_1}$ fixes $y$ we see
that $d_T(y, b_{i_1}a_{i_1}(y))=d_T(y, a_{i_1}(y))>k$.  
Let us look closer at what happens if $b_{i_1}\neq 1$.
The geodesic from $y$ to $a_{i_1}(y)$ goes through the vertex
$x(a_{i_1})$ and thus also through that vertex $y(b_{i_1})$.  Note that
the vertex $y(b_{i_1})$ is the vertex in
$[y,b_{i_1}a_{i_1}(y)]\cap[x,y]$ that is furthest away from $y$ and thus $d_T(y,\pr_{[x,y]}(h(y)))=d_T(y,y(b_{i_1}))<d_T(y,z)$.

Assuming the induction hypothesis above we write 
$h= b_{j_{l+1}}a_{i_{l+1}}b_{j_l}a_{i_l}\cdots b_{j_1}a_{i_1}c$ with all the $a_{i}$'s and $b_{j}$'s  occurring non-trivial except possibly $a_{i_1}$ and $b_{j_{l+1}}$.
Write $h'= b_{j_l}a_{i_l}\cdots b_{j_1}a_{i_1}c$ and note that $b_{j_{l}}\neq 1$.
The induction hypothesis says that $d_T(y,\pr_{[x,y]}(h'(y)))<d_T(y,z)$.
Observe that 
$$d(y,\pr_{[x,y]}(h'(y))+d(\pr_{[x,y]}(h'(y)),h'(y))=d(y,h'(y))>k.$$
The geodesic from $x$ to $h'(y)$ contains
$x(a_{i_{l+1}})$ and $\pr_{[x,y]}(a_{i_{l+1}}h'(y))=x(a_{i_{l+1}})$.  Clearly 
$d(\pr_{[x,y]}(a_{i_{l+1}}h'(y)),a_{i_{l+1}}h'(y))> d(\pr_{[x,y]}(h'(y)),h'(y))$ and $$d(y,\pr_{[x,y]}(a_{i_{l+1}}h'(y)))=d(y,x(a_{i_{l+1}}))> d(y,\pr_{[x,y]}(h'(y)).$$
    Hence
$d_T(y,a_{i_{l+1}}h'(y))> d(y,h'(y))>k$.  Therefore  $$d_Y(y,h(y))=d_T(y,b_{i_{l+1}}a_{i_{l+1}}h(y))>k.$$
Note that the geodesic from $y$ to
$a_{i_{k+1}}h(y)$ goes through the vertex $y(b_{i_{l+1}})$ and thus
$\pr_{[x,y]}(b_{i_{l+1}}a_{i_{l+1}}h(y))=y(b_{i_{l+1}})$ and therefore $d_T(y, \pr_{[x,y]}(h(y)))<d_T(y,z)$. 

\begin{theo}{\rm (Cf.~\cite[Theorem~2.5]{Smith2010})}\label{TSimon}
Let $G$ be a group acting on a tree $T$ with two orbits $V_1$ and $V_2$ on the vertex set.  Suppose that there are distinct vertices $x$ and $y$ in $V_1$ and  that on the path between $x$ and $y$ is a vertex $z$, distinct from both $x$ and $y$, such that $G_{x, z}=G_{y, z}$.   Then $G$ does not act primitively on $V_1$.
\end{theo}  

{\em Proof.}  Let $x, y$ and $z$ be as in the Theorem.  Suppose that $G$ acts primitively on $V_1$.  Then $G$ acts transitively on $V_1$ and, since $G_x$ is a maximal subgroup of $G$, then $\langle G_x, G_y\rangle=G$.  But now we have a contradiction with Lemma~\ref{LSimon} because $k=d_T(x,y)\geq 2$ and the orbit of $y$ under $\langle G_x, G_y\rangle=G$ would have to contain vertices in distance 2 from $y$ but by the Lemma that is impossible.  

\begin{coro}\label{CPrimitive}
Let $G$ be a group acting on a tree $T$ with two orbits $V_1$ and $V_2$ on the vertex set.   Suppose that $G$ acts primitively on $V_1$ and that the tree has infinite diameter.  Then the stabilizer $G_x$ of a vertex $x$ in $V_1$ is infinite.
\end{coro}

{\em Proof.}  Suppose the stabilizers of vertices in $V_1$ are finite.  Let $g$ be an element in $G$ that acts like a translation on $T$ and let $\{v_i\}_{i\in\Z}$ be the line $L$ that $g$ acts on by translation.  (Such an element exists by Lemma~\ref{LNoInvariant} part (iii).)  Assume that $v_0$ is in $V_1$.  Set $u_j=h^j(v_0)$.  Define $G(i)$ as the stabilizer of $u_i$ and $G(i,j)$ as $G_{u_i}\cap G_{u_j}$.  Note that $G(i,j)$ fixes all the vertices in the path between $u_i$ and $u_j$.  Hence $G(0)\supseteq G(0, 1)\supseteq G(0, 2)\supseteq\cdots$.  Because $G(0)$ is finite this sequence must eventually stop.  So there is a number $m$ such that $G(0, m)$ is equal to $G(0, j)$ for all $j\geq m$, i.e. the group $G(0, m)$ fixes all vertices $u_j$ with $j\geq m$.  Now $G(0, m)\supseteq G({-1}, m)\supseteq G({-2}, m)\supseteq\cdots$.  Since $G(0, m)$ is finite there is number $n\leq 0$ such that $G(n, m)=G(j, m)$ for all $j\leq n$.  The conclusion is that the group $G(n, m)$ fixes all the $u_i$'s and hence fixes all the vertices on the line $L$ and $G(n, m)$ is indeed equal to the pointwise stabilizer of the line $L$.  Note that $g^{m-n}G(n, m)g^{-(m-n)}=G(m,{2m-n})$.  We now set $x=u_n$, $z=u_m$ and $y=u_{2m-n}$, and see that $G_{x,z}=G_{y,z}$ and by Theorem~\ref{TSimon} it is now impossible that $G$ acts primitively on $V_1$.  We have reached a contradiction and therefore the assumption that the stabilizer of a vertex in $V_1$ is finite must be wrong.

\begin{coro}\label{CPrimitiveTree}
Let $G$ be a group acting primitively on a locally finite connected graph $X$ with infinitely many ends.  Then the  stabilizer $G_x$ of a vertex $x$ in $X$ is infinite and the group $G$ with the permutation topology is not discrete.
\end{coro}

{\em Proof.}   The action of $G$ on the tree $T_Y$ as described in Step 1 of the proof of Theorem~\ref{TPrimitiveSimple} satisfies the conditions in Corollary~\ref{CPrimitiveTree} with $V_1$ being the set of vertices in the tree $T$ that corresponds to the vertex set of $X$.  The stabilizer of a vertex $x$ in $T$ is equal to the stabilizer of the corresponding vertex in $X$ and the conclusion follows from Corollary~\ref{CPrimitive}.  Keeping in mind that the graph $X$ is locally finite we conclude that $G$ is not discrete.

\end{document}